\documentstyle[12pt,amstex,amssymb,amsfonts]{amsart}
\def\be#1{\begin{equation} \label{#1}}
\def\bs{\begin{split}}
\def\ba{\begin{align}}
\def\bas{\begin{align*}}
\def\R{{\mbox{\bf R}}}
\def\M{M}

\def\Z{{\mbox{\bf Z}}}
\def\eps{\varepsilon}

\def\xp{{x^\prime}}
\def\qp{{q^\prime}}
\def\rp{{r^\prime}}

\def\qpt{{\tilde q^\prime}}
\def\rpt{{\tilde r^\prime}}

\newcommand{\qtil}{{\tilde{q}}}

\newcommand{\rtil}{{\tilde{r}}}

\def\emph#1{{\it #1}}
\newenvironment{proof}{\noindent {\bf Proof} }{\endprf\par}
\def \endprf{\hfill  {\vrule height6pt width6pt depth0pt}\medskip}

\theoremstyle{plain}
  \newtheorem{theorem}[subsection]{Theorem}
  \newtheorem{proposition}[subsection]{Proposition}
  \newtheorem{lemma}[subsection]{Lemma}
  \newtheorem{corollary}[subsection]{Corollary}

\theoremstyle{remark}

\theoremstyle{definition}
  \newtheorem{definition}[subsection]{Definition}

\begin{document}

\pagestyle{empty}
\parindent = 0 pt
\parskip = 12 pt
\baselineskip = 18 pt
\vglue 1.8 cm
\centerline{\large SPHERICALLY AVERAGED ENDPOINT STRICHARTZ}
\centerline{\large ESTIMATES FOR THE TWO-DIMENSIONAL}
\centerline{\large SCHR\"ODINGER EQUATION\footnote{1991
\emph{Mathematics Subject Classification.} 35J10,42B25.}}
 
\centerline{}
\centerline{}
\centerline{Terence Tao}
\centerline{}
\centerline{Department of Mathematics,}
\centerline{UCLA, Los Angeles,}
\centerline{CA 90095-1555}
\centerline{\tt tao@@math.ucla.edu}
\centerline{}
\centerline{}
\centerline{}
\centerline{}
\begin{abstract}
The endpoint Strichartz estimates for the Schr\"odinger
equation is known to be false in two dimensions\cite{montgomery-smith}.  
However, if one averages the solution in $L^2$ in the angular variable,
we show that the homogeneous endpoint and the retarded
half-endpoint estimates hold, but the full retarded endpoint fails.
In particular, the original versions of these estimates hold for
radial data.
\end{abstract}

\section{Introduction}\label{intro}

Let $\Delta$ be the Laplacian on $\R^n$ for $n \geq 1$, so that $e^{it\Delta}$ is the
evolution operator corresponding to the free Schr\"odinger equation.

We consider the problem of obtaining bounds for this operator
in the mixed spacetime Lebesgue norms 

$$ \| F\|_{L^q_t L^r_x} = \int \| F(t,\cdot) \|_{L^r_x}^q\ dt)^{1/q}.$$

Such estimates
are commonly known as Strichartz estimates and have application to the study
of non-linear Schr\"odinger equations (see e.g. \cite{cwI}).
The following
Strichartz 
estimates are known (see \cite{ginebre:summaryschrodinger},\cite{tao:keel}):

\begin{definition} If $n$ is given, we say that the exponent
pair $(q,r)$ is \emph{admissible} if
$q,r \geq 2$, $(q,r,n) \neq (2, \infty, 2)$ and

$$
\frac{1}{q} + \frac{n}{2r} = \frac{n}{4}.
$$

\end{definition}

\begin{theorem}\cite{tao:keel}  If $n$ is given and 
$(q,r)$, $(\qtil,\rtil)$ are admissible, then we have the estimates

\begin{equation}\label{main-strichartz}
\| e^{it\Delta}f \|_{L^q_t L^r_x} \lesssim \|f\|_{L^2_x},
\end{equation}

\begin{equation}\label{dual}
\| \int e^{-is\Delta} F(s)\ ds \|_{L^2_x} \lesssim \|F\|_{L^\qp_t L^\rp_x},
\end{equation}

\begin{equation}\label{retarded}
\bigl\| \int_{s<t} e^{i(t-s)\Delta}F(s)\ ds\bigr\|_{L^q_t L^r_x}
\lesssim \| F\|_{L^\qpt_t L^\rpt_x}
\end{equation}

for all test functions $f$, $F$ on $\R^n$, $\R^{n+1}$ respectively.
\end{theorem}

The above conditions on $(q,r)$ are known to be necessary for the homogeneous
estimates \eqref{main-strichartz}
and \eqref{dual}, but it is not known what the necessary and sufficient conditions
are for the inhomogeneous retarded estimate \eqref{retarded}.  
For further discussion
we refer the reader to \cite{ginebre:summaryschrodinger}, \cite{montgomery-smith}, 
\cite{tao:keel}.  In the radial case for $n>2$ there is a further smoothing effect,
see \cite{vilela}.  

In this paper we investigate the ``forbidden endpoint'' $(q,r,n)=(2,\infty,2)$.
Accordingly we shall restrict ourselves to the two-dimensional case $n=2$
for the remainder of this paper.
With no further assumptions on $f$, $F$, the estimates
\eqref{main-strichartz}, \eqref{dual}, \eqref{retarded} are known to
be false even if $(\qtil, \rtil)$ are admissible, and even if the $L^\infty_x$
norm is replaced with the BMO norm; see \cite{montgomery-smith}.  The counterexamples
are non-radial and involve Brownian motion.

However, one can recover the endpoint estimate by averaging in $L^2$
over angular directions.  More precisely, let $L^\infty_r L^2_\theta$
denote the norm

$$ \| f \|_{L^\infty_r L^2_\omega} = \sup_{r>0} (\frac{1}{2\pi}\int_0^{2\pi} |f(re^{i\theta})|^2\ d\theta)^{1/2},$$

with the dual norm $L^1_r L^2_\theta$ defined similarly.  Then we have

\begin{theorem}  Let $(q,r,n)=(2,L^\infty_r L^2_\theta,2)$. 
Then
\eqref{main-strichartz} and \eqref{dual} hold, and the estimate 
\eqref{retarded} holds if $(\qtil,\rtil)$ is admissible.
\end{theorem}

For radial functions the $L^\infty_r L^2_\theta$ norm is just the $L^\infty$
norm, and so we have\footnote{After this paper was completed, we
learnt that this Corollary was independently proved by Atanas Stefanov.}

\begin{corollary}  Let $(q,r,n) = (2,\infty,2)$, and let $f$ and $F$ be radial.
Then
\eqref{main-strichartz} and \eqref{dual} hold.  The estimate \eqref{retarded}
holds if $(\qtil,\rtil)$ is admissible.
\end{corollary}

Finally, we present a very simple

\begin{proposition}  Let $(q,r,n) = (\qtil,\rtil,n) = (2,\infty,2)$.  Then
\eqref{retarded} can fail even if $F$ is radial.
\end{proposition}

This paper is organized as follows.  We first prove \eqref{main-strichartz}
for radial $f$ in Section \ref{homog}: the estimate \eqref{dual} follows
immediately by duality.  It turns out that the estimate reduces
easily to a maximal oscillatory integral estimate of the type discussed
in \cite{stein:kahane}, with a minor complication arising from
the behaviour of the Bessel function 

\begin{equation}\label{bessel-def}
J_n(x) = \frac{1}{2\pi} \int_0^{2\pi} e^{ix \cos\theta} e^{in\theta}\ d\theta
\end{equation}

for $x$ close to $n$.

We then turn to the positive results for \eqref{retarded} in Section
\ref{ret}.  Fortunately 
we shall be able to obtain these results as an automatic consequence
of the homogeneous estimate, by a very general argument
of Christ and Kiselev \cite{christ:kiselev}.

Finally, we discuss the negative results in Section \ref{negative}.

We remark that analogous results hold for the wave equation (with
$n=3$ playing the role of $n=2$) but are proved differently.  See 
\cite{tao:keel},
\cite{klainerman:nulllocal},
\cite{montgomery-smith}. 

\section{The homogeneous estimate}\label{homog}

In this paper $C$ denotes an absolute positive constant which may vary
from line to line, and we use the notation $A \lesssim B$ as shorthand
for $A \leq CB$.

In this section we prove \eqref{main-strichartz} for $(q,r,n) =
(L^2,L^\infty_r L^2_\theta,2)$, which implies \eqref{dual}
by duality.

We will always make the a priori assumption that all functions
are in the Schwarz class, and any singular integrals will be evaluated
in the principal value sense.

Since $e^{it\Delta}$ commutes with rotations, and our norms are
$L^2$ in the angular variable, we may use standard orthogonality 
arguments to
reduce to the case when $f$ is
given by a single spherical harmonic, i.e.

$$  f(re^{i\theta}) = f_n(r) e^{in\theta}$$

for some $n \in \Z$ and some function $f_n(r)$.  Our task is then to 
prove \eqref{main-strichartz} with
a bound independent of $n$.  

Fix $n$; we may assume that $n \geq 0$.  From the explicit formula for the 
fundamental solution

\be{explicit}
 e^{it\Delta} f(x) = \frac{C}{t} \int e^{i|x-y|^2/4t} f(y)\ dy
\end{equation}

and a change to polar co-ordinates, we obtain

$$
e^{it\Delta} f(r e^{i\theta}) = \frac{C}{t} 
\int e^{i|re^{i\theta} - Re^{i\phi}|^2/4t} f(Re^{i\phi})\ d\phi\ RdR,$$

which simplifies to

$$ e^{it\Delta} f(r e^{i\theta}) = \frac{C}{t} e^{i r^2/4t} \int_0^\infty \int_0^{2\pi}
e^{i rR\cos(\theta-\phi)/2t} e^{iR^2/4t} f_n(R)e^{in\phi}\ d\phi\ RdR.$$

Making a change of variables $\alpha = \theta-\phi$ and taking
absolute values, we obtain

$$ |e^{it\Delta} f(r e^{i\theta})| =
\frac{C}{|t|} \bigl|\int_0^\infty (\int_0
^{2\pi}
e^{i rR\cos(\alpha)/2t} e^{in\alpha} d\alpha) e^{iR^2/4t} f_n(R)\ RdR
\bigr|.$$

By \eqref{bessel-def}. the inner integral is $2\pi J_n(\frac{rR}{2t})$.
Thus \eqref{main-strichartz} can be rewritten as

$$ (\int \sup_{r \geq 0} |\int_0^\infty J_n(\frac{rR}{2t}) e^{iR^2/4t} 
f_n(R)\ RdR|^2  \ \frac{dt}{t^2})^{1/2}
\lesssim (\int_0^\infty |f_n(R)|^2\ RdR)^{1/2}$$

From our a priori assumptions we may replace $r \geq 0$ by $r>0$ in the supremum.

Write $\xi = R^2$, and $g(\xi) = f_n(R)$.  Also write $x = 1/(8 \pi t)$
and $\lambda = r/(2|t|)$.
After a change of variables, the above estimate becomes

$$
(\int \sup_{\lambda > 0} |\int_0^\infty J_n(\lambda \xi^{1/2})
e^{2 \pi i x \xi} g(\xi)\ d\xi|^2\ dx)^{1/2} \lesssim
(\int_0^\infty |g(\xi)|^2\ d\xi)^{1/2}.$$

Clearly this estimate will be implied by

$$
(\int \sup_{\lambda > 0} |\int J_n(\lambda |\xi|^{1/2})
e^{2 \pi i x \xi} g(\xi)\ d\xi|^2\ dx)^{1/2} \lesssim
(\int |g(\xi)|^2\ d\xi)^{1/2}$$

where the integrations now range over all of $\R$.

Let $G$ be the Fourier transform of $g$.  By Plancherel's theorem,
the above estimate is equivalent to

\be{t-targ}
 \| \sup_{\lambda > 0} |T_\lambda G| \|_2 \lesssim \|G\|_2
\end{equation}

where $T_\lambda$ is the multiplier defined by

$$ \widehat{T_\lambda G}(\xi) = J_n(\lambda |\xi|^{1/2}) \hat{G}(\xi).$$

We partition the Bessel function $J_n$ smoothly as

$$ J_n(r) = m_0(r) + m_1(r) + \sum_{2^j \gg n} m_j(r)$$

where $m_0$, $m_1$, and $m_j$ are supported on $|r| \ll n$, $|r| \sim n$,
and $|r| \sim 2^j \gg n$ respectively.  We similarly decompose $T_\lambda$ as

$$ T_\lambda = T^0_\lambda + T^1_\lambda + \sum_{2^j \gg n} T^j_\lambda.$$

Finally, for $j=0,1$ or $2^j \gg n$ we let $K^j_\lambda$ be
the convolution kernel of the operator $T^j_\lambda$; note that

$$ K^j_\lambda(x) = \int e^{2\pi i x \xi} m_j(\lambda |\xi|^{1/2})\ dx.$$

When $|r| \ll n$, the phase in \eqref{bessel-def} is non-stationary.
From this one can easily obtain the bounds

$$ \| m_0 \|_{C^k} \lesssim n^{-N}$$

for any $N, k > 0$.  Since $m$ is compactly supported, we thus have

$$ |K^0_\lambda(x)| \lesssim n^{-N} \lambda^{-2} (1 + \lambda^{-2} x)^{-k}$$

for any $N, k > 0$.  This in turn implies that

$$ |T^0_\lambda G(x)| \lesssim n^{-N} \M G(x)$$

where $\M$ is the Hardy-Littlewood maximal function.  Thus the contribution
of $T^0_\lambda$ to \eqref{t-targ} is acceptable.

We now turn to the contribution of $T^1_\lambda$.  We will not attempt
to estimate this contribution efficiently, and rely instead on very
crude tools.  We begin with
the Sobolev embedding $H^1(\R) \subset L^\infty(\R)$, which we write as

$$ \sup_y |f(y)| \lesssim (\int |f(y)|^2 + |f^\prime(y)|^2\ dy)^{1/2}.$$

We apply the change of variables $\lambda = e^{y/n}$ to obtain

$$ \sup_\lambda |g(\lambda)| \lesssim
(\int (n|g(\lambda)|^2 + \frac{1}{n}|\lambda g^\prime(\lambda)|^2) \frac{d\lambda}{\lambda})^{1/2}.$$

Applying this to $g(\lambda) = T^1_\lambda G(x)$ and taking $L^2$ norms of
both sides,
we obtain

$$
\| \sup_\lambda |T^1_\lambda G| \|_2 \lesssim
(\int (n \|T^1_\lambda G\|_2^2 + \frac{1}{n} \|\lambda \frac{\partial}{\partial\lambda} 
T^1_{\lambda} G\|_2^2)\ \frac{d\lambda}{\lambda})^{1/2}.$$

Applying Plancherel's theorem, we see that we will be done provided that

$$ 
\int (n |m_1(\lambda |\xi|^{1/2})|^2 + \frac{1}{n} |\lambda \frac{\partial}{\partial \lambda} m_1(\lambda |\xi|^{1/2})|^2) \frac{d\lambda}{\lambda} \lesssim 
1$$ 

uniformly in $\xi$.
By rescaling $\lambda$ by $|\xi|^{1/2}$ 
we may assume $\xi = 1$; from the support of $m_1$
we can thus restrict the integration to the region $\lambda \sim n$.  
It thus suffices to show that

$$
\int_{\lambda \sim n}
|m_1(\lambda)|^2 + |m_1^\prime(\lambda)|^2 \lesssim 1.$$

However, from \eqref{bessel-def}, the definition of $m_1$,
and Van der Corput's lemma (See e.g. \cite{stein:large}) we have
the estimates
\bas
|m_1(\lambda)| &\lesssim n^{-1/3} (1 + n^{-1/3} |\lambda - n|)^{-1/4}\\
|m_1^\prime(\lambda)| &\lesssim n^{-1/2},
\end{align*}
and the claim follows.

Finally, we consider the contribution of the $T_\lambda^j$ to \eqref{t-targ}.
We will show

$$ \| \sup_{\lambda > 0} |T_\lambda^j G| \|_2 \lesssim 2^{-j/4}
\|G\|_2$$

uniformly for $j$, $n$ such that $2^j \gg n$; this clearly implies
that the contribution of the $T_\lambda^j$ is acceptable.

Fix $j$.  It suffices to show that

$$ \| T_{\lambda(x)}^j G\|_{L^2_x} \lesssim 2^{-j/4} \|G\|_2$$

for an arbitrary positive function $\lambda(x)$ which we now fix.  We write this
as 

$$ \| \int G(y) K_{\lambda(x)}^j(x-y)\ dy\|_{L^2_x} 
\lesssim 2^{-j/4} \|G\|_2$$

By the $TT^*$ method, it suffices to show that

\be{targ}
 \| \int(\int K_{\lambda(x)}^j(x-y) \overline{K_{\lambda(\xp)}^j(\xp-y)}\
dy) F(\xp)\ d\xp \|_2 \lesssim 2^{-j/2} \|F\|_{L^2_\xp}
\end{equation}

for all test functions $F$.
This will follow from the estimate

\begin{lemma}  For any $a,b > 0$, $x, \xp \in \R$ we have

$$
|\int K_{a}^j(x-y) \overline{K_{b}^j(\xp-y)}\ dy|
\lesssim \Phi_{j,a}(|x-\xp|)
$$

where $\Phi_{j,a}$ is a radial decreasing non-negative function with

$$ \sup_a \| \Phi_{j,a} \|_1 \lesssim 2^{-j/2}.$$

\end{lemma}

Indeed, from this lemma we may bound the left-hand side of
\eqref{targ} pointwise by $C 2^{-j/2} \M F(x)$.

\begin{proof}
By Parseval's theorem and the definition of $K_\lambda^j$,  
the left-hand side is

\be{mj-int} C |\int m_j(a |\xi|^{1/2}) \overline{m_j(b |\xi|^{1/2})}
e^{2\pi i (x-\xp)\xi}\ d\xi|.
\end{equation}

On the other hand, from the standard asymptotics of $J_n$ (see e.g. \cite{stein:large}) we have

$$ m_j(\xi) = \sum_\pm 2^{-j/2} e^{\pm i|\xi|} \psi^\pm_j(2^{-j} \xi)$$

where $\psi^\pm_j(\xi)$ is a bump function on $|\xi| \sim 1$ uniformly
in $j$, $n$.  We can therefore rewrite \eqref{mj-int} as a finite
number of expressions of the form

$$ C 2^{-j} | 
\int e^{i(\pm a \pm b) |\xi|^{1/2}} e^{2\pi i (x-\xp)\xi} 
\psi^\pm_j(2^{-j} a |\xi|^{1/2})
\overline{\psi^\pm_j(2^{-j} b|\xi|^{1/2})}\ d\xi|
$$
where the $\pm$ signs need not agree.

It suffices to estimate the $\xi > 0$ 
portion of the integral.  From the change of variables 
$\xi = 2^{2j} a^{-2} t^2$, this
becomes

$$
C 2^j a^{-1} | \int e^{2\pi i (\frac{2^j(\pm a \pm b)}{2\pi a}t + \frac{2^{2j} (x-\xp)}{a^2} 
t^2)} 
\psi_\pm(t) \psi_\pm(\frac{b}{a} t)\ tdt|.
$$

We will majorize this expression by $\Phi_{j,a}(|x-\xp|)$, where

$$ \Phi_{j,a}(r) = C \min(r^{-1/2}, 2^j a^{-1}, 2^j a^{-1} (2^j a^{-1} r)^{-10});$$

it is easy to verify that $\Phi$ satisfies the desired properties.

The first bound of $C r^{-1/2}$ follows from Van der Corput's 
lemma (see e.g. \cite{stein:large}).  The second bound of $C 2^j a^{-1}$ simply
follows from replacing everything by absolute values.  To show
the third bound, we may assume from the second bound that
$|x-\xp| \gg 2^{-j} a$.  But then the phase
$\frac{2^j(\pm a \pm b)}{2\pi a}t + \frac{2^{2j} (x-\xp)}{a^2}$ has
a derivative of magnitude at least $2^j a^{-1} r$ 
on the support of $\psi_\pm$, and the bound follows from
non-stationary phase.
\end{proof}

One can improve the estimate on $T^1_\lambda$ by incorporating
techniques from the treatment of $T^0_\lambda$ and $T^j_\lambda$.  This
will eventually yield a gain of $n^{-\eps}$ to \eqref{t-targ} for
some $n > 0$.
This translates to a gain of angular regularity, so that
we may replace $L^2_\theta$ by an angular Sobolev space $H^\eps_\theta$.  
By Sobolev embedding this implies that the $L^2_\theta$ can be improved to
$L^p_\theta$ for some $p>2$.  A possibly related
smoothing effect in higher dimensions was observed in \cite{vilela}.

It is not clear what the best value of $p$ is.  However 
the negative results in \cite{montgomery-smith} show that this cannot be
improved to $p=\infty$ or to $p=\rm{BMO}$.

\section{The inhomogeneous estimate}\label{ret}

We now prove \eqref{retarded} when $(q,r,n) = (2, L^\infty_r L^2_\theta,2)$
and $(\qtil, \rtil)$ is admissible.  We first observe that if the
restriction $s<t$ were somehow removed from the integral, the left-hand side
of \eqref{retarded} would factor as

$$ \int e^{i(t-s)\Delta}F(s)\ ds = e^{it\Delta} ( \int e^{-is\Delta} F(s)\ ds),$$

and the claim would then follow by 
combining \eqref{main-strichartz} and \eqref{dual}.

To finish the proof we need to reinstate the restriction $s<t$.  This
can in fact be done very general circumstances, as observed by
Christ and Kiselev \cite{christ:kiselev}, \cite{christ:kiselev2}.  
More precisely, we have

\begin{lemma}\cite{christ:kiselev}  Let 

$$Tf(t) = \int_\R K(t,s) f(s)\ ds$$ 

be a linear transformation which maps
$L^p(\R)$ to $L^q(\R)$ for some $1 < p < q < \infty$.  Then the map

$$ \tilde T f(t) = \int_{s < t} K(t,s) f(s)\ ds$$

also maps $L^p(\R)$ to $L^q(\R)$. 
\end{lemma}

For our purposes we need the trivial observation
that the argument below extends to the case
when $K$ takes values in $B(X,Y)$, the space of bounded mappings
from one Banach space to another.

\begin{proof}  
We will prove the claim for smooth $f$ only, to avoid technical problems.
We normalize so that $\|f\|_p = 1$.
Define the function $F(t)$ by $F(t) = \int_{s < t} |f(s)|^p\ ds$.  This
map $F$ is an order-preserving bijection from $\R$ to $[0,1]$.

Partition the interval $[0,1]$ into dyadic intervals in the usual manner.
We define a relationship $I \sim J$ on dyadic intervals as follows:
$I \sim J$ if and only if $I$ and $J$ are the same size, are adjacent, 
and the elements of $I$ are strictly
less than the elements of $J$.  It is easy to verify that for almost every
$x<y$ there is a unique pair $I$, $J$ such that $x \in I$, $y \in J$, and
$I \sim J$.  Applying this with $x=F(s)$, $y=F(t)$, we obtain

$$
\int_{s<t} \ ds =
\int_{F(s)<F(t)} \ ds =
 \sum_{I,J: I \sim J} \chi_{F^{-1}(J)}(t) 
\int_{F^{-1}(I)}\ ds.$$

We thus have

$$  \tilde T f = \sum_{I,J: I \sim J} \chi_{F^{-1}(J)} T (\chi_{F^{-1}(I)} f).$$

We need to show that $\|\tilde T f\|_q \lesssim 1$.  It suffices to prove that

\be{t}
  \|\sum_{I,J: I \sim J, l(I) = 2^{-j}} \chi_{F^{-1}(J)} T (\chi_{F^{-1}(I)} f)\|_q \lesssim 2^{-\eps j}
\end{equation}

uniformly in $j \geq 0$ for some $\eps > 0$, where $l(I)$ denotes the
sidelength of $I$.

Fix $j$.  Since for each $I$ there are at most two $J$, and the 
functions $\chi_{F^{-1}(J)}$ have essentially disjoint support, we can
estimate the left-hand side of \eqref{t} by

$$ (\sum_{I: l(I) = 2^{-j}} \| T (\chi_{F^{-1}(I)} f)\|_q^q)^{1/q}.$$

By the assumption on $T$, this is bounded by

$$ (\sum_{I: l(I) = 2^{-j}} \| \chi_{F^{-1}(I)} f\|_p^q)^{1/q}.$$

But by construction $\| \chi_{F^{-1}(I)} f\|_p = 2^{-j/p}$, hence this
sum is just

$$ 2^{-j(\frac{1}{p}-\frac{1}{q})},$$

and the claim follows from the hypothesis $p<q$.
\end{proof}

The requirement $p < q$ is necessary, as can be seen by considering
the Hilbert transform.  The lemma also holds in the ranges
$1 = p \leq q \leq \infty$ and $1 \leq p \leq q = \infty$, but for
more trivial reasons.  We remark that a stronger maximal version of this 
lemma appears in \cite{christ:kiselev}.

\section{Negative results}\label{negative}

We now show why \eqref{retarded} fails when $(q,r,n)=(2,\infty,2)$ and
$(\qtil,\rtil)$ is not admissible, even when $F$ is radial.

From dimensional analysis (recalling that time has twice the dimensionality
of space for the purposes of the Schr\"odinger equation) 
we obtain the necessary condition for \eqref{retarded}

$$ \frac{2}{q} + \frac{2}{r} + 2 = \frac{2}{\qtil'} + \frac{2}{\rtil'}.$$

Thus we must have 

$$\frac{1}{\qtil}+\frac{1}{\rtil} = \frac{1}{2}.$$

Therefore the only case left to consider is the double forbidden endpoint

$$ (q,r,n) = (\qtil,\rtil,n) = (2,\infty,2).$$

By a limiting argument we may assume that $F$ is a measure on the time axis
$x=0$:

$$ F(x,s) = g(s) \delta(x).$$

Since $G(0) \leq \|G\|_\infty$ for any $G$, 
it suffices to disprove the estimate

$$
\bigl\| \int_{s<t} [e^{i(t-s)\Delta}F(s)](0)\ ds\bigr\|_{L^2_t}
\lesssim \| g\|_{L^2_s}.
$$

By \eqref{explicit}, this is

$$ 
\bigl\| \int_{s<t} \frac{1}{s-t} g(s)\ ds \|_{L^2_t} \lesssim \|g\|_{L^2_s},$$

which is clearly false (e.g. take $g = \chi_{[0,1]}$).
 
It is easy to modify this argument to show that the estimate continues to fail
if the $L^\infty$ or $L^1$ norms are replaced by BMO or $H^1$ norms, or if
some frequency restriction or smoothness condition is placed on $F$.

\centerline{}
\centerline{\bf ACKNOWLEDGEMENTS}

This paper arose out of discussions with Tony Carbery
and Christoph Thiele.  The author thanks Elias Stein
and Steve Wainger
for pointing out the argument in \cite{stein:kahane}, and 
Michael Christ for the argument in \cite{christ:kiselev},
both of which greatly simplified this paper.  The author 
also thanks Tony
Carbery, Ana Vargas, and Sergiu Klainerman for their
hospitality while this work was performed.  The author
is partially supported
by NSF grant DMS-9706764.

\def\refname{\centerline{REFERENCES}}
\bibliographystyle{plain}

\end{document}